\input amstex
\documentstyle{amsppt}
\input epsf
\input psfrag

\define\RR{\Bbb R}
\define\BB{\Bbb B}
\define\CC{\Bbb C}
\define\SS{\Bbb S}
\define\FF{\Bbb F}

\hyphenation{Liou-ville}

\NoBlackBoxes

\topmatter
\title Complex Axis and de Medeiros\rq \ Campo Vetorial
\endtitle
\rightheadtext {}
\author Jon A. Sjogren
\endauthor
\affil Towson University
\endaffil
\address Towson, Maryland
\endaddress
\date 1 October 2018
\enddate
\dedicatory To the Memory of Robert Todd Gregory (1920--1984)\\ Pioneer of Exact Computation \enddedicatory
\endtopmatter
\document

\abovedisplayskip=13pt
\belowdisplayskip=13pt

\head Abstract\endhead

In linear algebra a fundamental question arises: does any endomorphism $\frak M: \CC^n \to \CC^n$ have an {\it axis} (complex vector of invariance)? The {\it real} case $A: \RR^n \to \RR^n$ is well-understood when $n$ is an {\it odd} integer. A real root of the characteristic polynomial $\sigma_A (t)$ exists, which yields an eigen-value. Since $\sigma_A$ amounts to a non-linear construction, can one settle this {\it odd real} and also the {\it complex} case in a more geometric manner?

An elegant solution was put forward by A. von Sohsten de Medeiros \cite{de Medeiros}. The proof of de Medeiros shows the existence of a complex axis for non-singular $\frak M$ by applying the Lefschetz fixed-point theorem to a related continuous mapping $\tilde{\frak M}: \CC P^{n-1}\to \CC P^{n-1}$. An eigen-vector for a singular matrix is found through standard row reduction. The induced mapping $\tilde{\frak M}$ is well-defined so long as $\frak M$ gives a $\CC$-isomorphism (invertible).

We wish to work with the main insights of de Medeiros' approach, but aim for a dual implementation. Instead of inducing a family of continuous mappings, we define a vector field $v_{\frak M}$ on $\CC P^{n-1}$. When $\frak M$ affords no eigen-vector, $v_{\frak M}$ is a never-vanishing field. But this cannot be since the Euler characteristic $\chi (\CC P^{n-1}) = n > 0$. We try to avoid any calculation based on the Lefschetz number or Euler characteristic. We stop short of the Poincar\'e--Hopf theorem, and avoid stochastic (or ``measure'') concepts common in differential topology such as the theorem of Sard, or general-position arguments. Instead we use Hopf's Lemma on the invariance of total index of a non-degenerate vector field (a field whose zeros are ``non-degenerate''). At this point, for comparison, {\it any} non-degenerate vector field on $\CC P^{n-1}$ can be selected, one that is easy to work with. Without much computation, one may confidently assert that the Milnor--Hopf vector field has total index $> 0$ (actually it equals $n$ as expected). This shows how our original matrix $\frak M$ violates the Hopf Lemma.

Then of course the solution to any monic polynomial equation $p(z) = 0$ of degree $n$, follows from the correspondence between a ``characteristic'' (or ``secular'') polynomial and its ``companion matrix''. We must exercise care that no part of the argument uses results derived from the ``Fundamental Theorem of Algebra''.

\head Introduction\endhead

Let $B: \CC^n \to \CC^n$ denote a linear endomorphism, which is also a surjection (hence an isomorphism). The question arises whether $B$ necessarily possesses an invariant vector $v \in \CC^n$, generating a complex line of invariance (axis), and can this fact be shown without appealing to the Fundamental Theorem of Algebra? This last means of course that any monic polynomial over the complexes $\CC$, of degree $n \geq 1$, can be decomposed completely into linear factors. To see that this point of view is not a frivolous one, consider the case of $n = 2m+1$, odd. Hence the method is to consider linear representations on a real vector space of matrices. This approach was proposed by Derksen and simplified in \cite{Conrad} to obtain the Axis Theorem for all $n \geq 1$, not only for odd degrees.

This proof and another one from \cite{Sjogren,\,Endo} are summarized in Section 1.  Now A. von Sohsten de Medeiros gives the same result using the rudiments of Algebraic Topology on projective spaces. In \cite{de Medeiros}, the author shows how an induced mapping $\tilde{B}: \CC P^{n-1}\to \CC P^{n-1}$ must be homotopic to the identity of this projective space, else a singular matrix arises whose null-vector fills the bill (as axis of invariance). Hence the Lefschetz number $\Cal L(\tilde{B})$ would have to equal $\chi (\CC P^{n-1}) = n \neq  0$. Therefore $\tilde{B}$ affords a fixed point $\xi \in \CC P^{n-1}$, hence the axis (complex line) $\hat{\xi}\subset \CC^n$ is held invariant by the original transformation $B$.

In the spirit of the de Medeiros methodology, will consider a dual version. Let us enumerate some advantages of this modification, and note a potential drawback.

Firstly, there is no distinction between singular and non-singular cases: we treat an endomorphism $B$ of lesser rank the same as one of full rank.

Secondly, there is no ``exceptional case'' as when $t B + (1-t)I$ turns out to be an $n \times n$ matrix of deficient rank (yielding the desired axis immediately).

Thirdly, we do not induce a mapping such as $\tilde{B}$ on the quotient space $\CC P^{n-1} \simeq \CC^n/\CC^{*}$ at all. Thus it is not necessary to be concerned with questions of general topology involving the continuity of a ``homotopy'' on this quotient space (defined under an ``identification mapping''). A typical example of this issue is discussed in \cite{Maunder}, p\. 19.

Finally, our original motivation for reconsidering the problem was that the Lefschetz--Euler principle of fixed points of a self-mapping does not seem the most natural way to look at projective space, geometrically. The proof of the Lefschetz FPT best known to the Applied Topologist starts with a simplicial decomposition of some space $X$ of interest, together with an approximation of the given self-mapping. Now the Hopf\footnote{We mean Heinz Hopf, 19 Nov 1894--3 June 1971.} Trace Formula asserts that  the alternating sum of traces over simplices of various dimensions, equals the alternating sum of the induced (``rational'') homology traces. Due to the fluctuating $\pm$ signs, terms in the sum of the {\it chain level} cancel, namely terms corresponding to the $\mbox{trace}(f)$ on $(k-1)$-dimensional boundary chains, and $\mbox{trace}(f)$ on $k$-dimensional chains, modulo cycles. But  the usual way to depict $\CC P^{n-1}$ geometrically is as a union of even-dimensional cells ($k$-cubes).
$$e_0 \cup e_2 \cup \cdots \cup e_{2(n-1)}$$
where $e_{2(j+1)}$ is attached to $e_{2j}$ via its bounding $\SS^{2j+1}$, called without ennui, the ``Hopf'' mapping.

In computing the trace on homology say at $e_2$, we observe {\it no} ``boundary chains'', of odd dimension, to cancel! It appears that an effective simplicial approximation of a self-map $f: X \to X$ requires the presence of simplices of odd dimension. To save on reading time, we point to the relevant one-line statements in \cite{Maunder}, \linebreak p\. 150 and \cite{Spanier}, p\. 172, Theorem 11.

By contrast, our present exposition invokes several concepts of geometry, including ``tubular neighborhood'' and ``degree'' of a vector field.  Instead of appealing to the more technical Poincar\'e--Hopf theorem, we quote {\it Hopf's Lemma}, which roughly computes the total index of a vector field on a manifold as the Brouwer--Heinz degree of the Gau{\ss} (``normal'') mapping of the boundary sub-manifold.  Thus this total index, where it is a finite sum, does not even depend upon the particular vector field which is being examined.

\head 1. The Axis in $\RR^{2m+1}$ and $\CC^{2m+1}$\endhead

K. Conrad provides an interpretation to the ``linear algebra proof'', due to H. Derksen, that any endomorphism of an odd-dimensional complex vector space affords an eigen-vector. In this treatment, see \cite{Conrad}, certain {\it real} representations are constructed. As a starting point, the eigen-vector or ``axis'' theorem is stated verbatim over the field $\RR$ of real numbers. This result amounts essentially to the ``cueball'' or ``Hedgehog'' principle, which is well-expressed geometrically by Milnor in \cite{CUE}, owing also to work in foliation theory of Dan Asimov.

A geometric attack on ``Derksen's Lemma'' (so-called since it is the base instance of an induction leading to the full Complex Axis Theorem) is given in \cite{Sjogren,\,Endo} and will be recapitulated.

Let $A: \RR^n \to \RR^n$ be a linear transformation where $n$ is odd. If $A$ is ``singular'', a sequence of row operations will put the corresponding square matrix $\hat{A}$ into an echelon form exhibiting a row of zeros. This leads to the construction of a ``null-vector'' which can be expressed by means of the original basis of $\RR^n = \RR^{2m+1}$.

The remaining case is where $\hat{A}$ has full rank $n$. Considered as a restricted mapping $\SS^{2m}\to \RR^n$, a scaled continuous mapping
$\bold a : \SS^{2m} \to \SS^{2m}$ arises according to
$$\bold a (s) =\frac{A(s)}{\|A(s)\|}\quad .$$

Here $\SS^{2m} = \left\{s \in\RR^n \,\, | \,\, \|s\| =1\right\}$. If $A$ has no real eigenvector, there is no $y \in \SS^{2m}$ such that $\bold a (y) = \pm y$. Considering the $\bold a (y)$ as a vector field in $\RR^n$ (each vector based at $y \in \SS^{2m}$), none of these vectors is collinear to $y \in \RR^n$ itself. Hence there will be a unique continuous projection of $\sigma (y)$ to a vector of unit norm in $T_y$, the tangent space to $\SS^{2m}$ at $y$. According to the Milnor-Asimov ``Hedgehog'' result (previously known as the ``Poincar\'e--Brouwer Theorem'', see \cite{Dugundji}), a vector field of this nature cannot exist. Thus any real square matrix of odd order affords a real eigen-vector, with real eigen-value. See Figure 1.

\bigskip
\centerline{\epsfbox{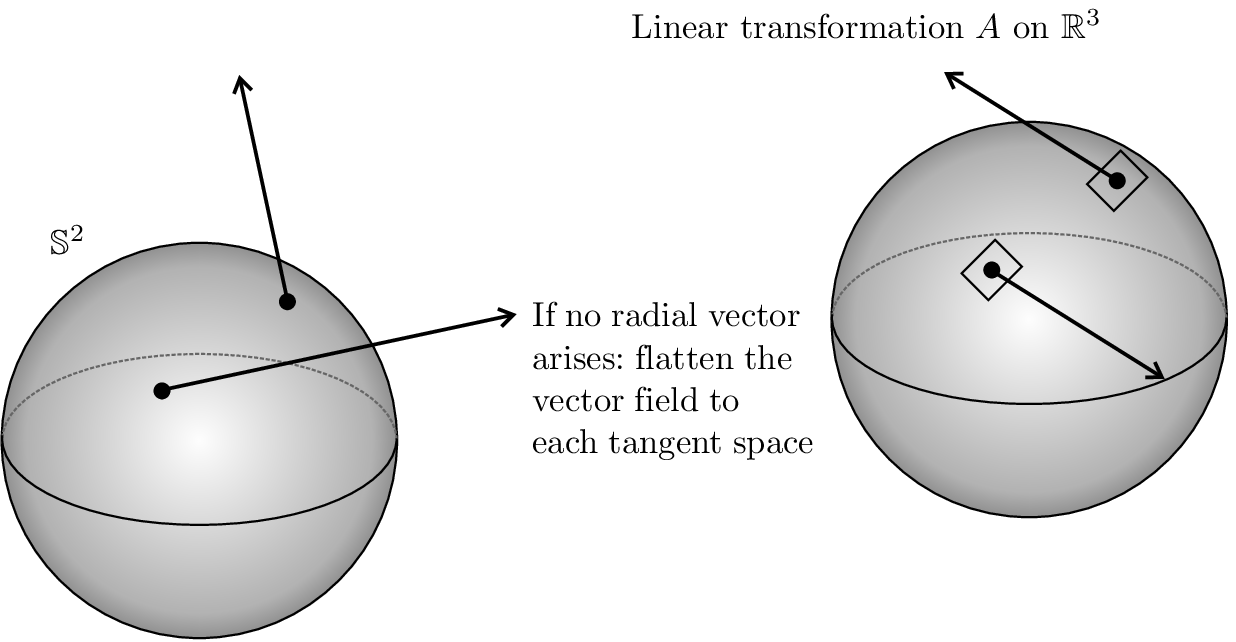}}

\medskip
\centerline{Figure 1. ``Hedgehog'' Violation}

\bigskip

Next we give a synopsis of the Conrad--Derksen proof for an ``odd complex'' endomorphism. Given linear $A: \CC^n \to \CC^n$, take $\hat{A}$ to be its $n\times n$ matrix in standard coordinates.
Let $B \in H_n$, the set of $n\times n$ Hermitian matrices, considered as an $n^2$-dimensional vector space over $\RR$. Let $L, K$ be the real operators defined by
$$L(B) =\frac{1}{2}\bigl(\hat{A}B + B\hat{A}^*\bigr), \quad K(B) = \frac{1}{2i}\bigl(\hat{A}B-B\hat{A}^*\bigr)\,\quad.$$
Here $\hat{A}^*$ is the conjugate transpose of $\hat{A}$. It is seen that those two real operators commute: for all $B \in H_n$, $L(K(B)) = K(L(B))$.

A Proposition from \cite{Derksen} states: Let $\FF$ be a field (as in algebraic structure). If every endomorphism $T$ on an {\it odd-dimensional} vector space $V$ over $\FF$ has a (non-zero) eigen-vector, then every {\it pair} $K,L$ of {\it commuting} endomorphisms on any {\it odd-dimensional} vector space $W$ has a {\it common} eigen-vector.

That is, there exists $w \in W$ and $\lambda_1, \lambda_2 \in \FF$ with $K(w) = \lambda_1\, w$, $L(w) = \lambda_2\, w$. In our case, focusing in $\FF=\RR$ and the vector space $H_n$ of Hermitian matrices, we conclude that given $\hat{A}$, there exists $B \in H_n$ so that for some $\lambda ,\, \kappa \in \RR$ there hold
 $$L(B) = \lambda B,\quad K(B) = \kappa B,$$
and $\hat{A}B = \left(\lambda + i\, \kappa\right) B$ (Exercise). Thus any non-zero column vector of $B \ne [\bold 0]$ can be taken as the eigenvector $b \in \CC^n$ that we are seeking, given
$\hat{A}\in M_n (\CC )$. Derksen's Proposition applies since we showed (geometrically) that any $T \in M_n (\RR)$ for $n = 2m+1$ must have a real eigen-vector.

 What we have said now proves a {\it complex} version of Milnor's ``Hedgehog'' factum. Indeed, a complex square matrix of odd order called $T$, that has no non-vanishing eigen-vector, leads to three real square matrices $I, D, E$ of ``twice the order'' $(2n)$, that together span {\it only} non-singular matrices, plus the all-zeros matrix. Such a circumstance is however impossible, as shown by work in homotopy theory done at about the same time by B. Eckmann in Switzerland and G. Whitehead at Chicago. Their result was that a $4m+1$-sphere has at most {\it one} linearly independent tangent vector field. The modern proof of this uses the cohomology ring of $\RR P^{2n-1}$. See \cite{Eckmann} and \cite{Whitehead} for greater detail.

 The study of linear spaces of matrices with bounded rank (including the square $\sigma$-matrix) is a well-developed theory, see \cite{Meshulam} and references therein. For our purposes, the important conclusion is

\noindent
{\bf Proposition 1}\quad  Any three real (square) matrices of order ``twice an odd'', say $A, B, C$ can be combined linearly to form a {\it singular} (non-invertible) matrix, using (the three) coefficients which are not all zero.

\noindent
{\it Proof} \quad Suppose firstly that no such singular linear combination exists. If $q$ denote the matrix order, then this {\it theory} that we refer to points to a line-bundle equality (isomorphism) over the two-dimensional manifold $\RR P^2$, namely
$$
\gamma^q \simeq \epsilon^q\qquad \mbox{(Whitney sum)}\tag{$\Delta$}\,\,\quad ,
$$
where $\epsilon$ is the trivial line bundle and $\gamma$ is the ``tautological'' line bundle.

A Lemma simplifies these equations: over $\RR P^2$ there holds $\gamma \oplus \gamma \oplus \gamma \oplus \gamma \simeq \epsilon \oplus \epsilon \oplus \epsilon \oplus \epsilon $. This is proved by constructing four independent sections of $\gamma^4$, based on the ``quaternion multiplication table''.

Because the base space has such a low dimension, considerable cancellation in $(\Delta)$ can now be effected. In fact, we are only required to show that
$$\gamma \oplus \gamma \oplus \epsilon \simeq \epsilon^3$$
is {\it not} possible see \cite{Husem\"oller}. This is a result sufficient to conclude that the vector bundle $\gamma^2$ is not {\it stably trivial}.

The geometric content of this final observation is just the {\it classical} Borsuk-Ulam theorem, which says that any continuous flattening of $\SS^2$ to $\RR^2$ sends some pair of antipodal points to the same image point. For an elementary treatment, see \cite{Sieradski}. Supposing that we had sections that realize a trivialization of $\gamma \oplus\gamma\oplus\epsilon$, we arrive at two self mappings $p_1,\, p_2: \SS^2 \to \SS^2$ which are {\it odd maps} (antipode-preserving), satisfying both $p_1 \sim p_2$ (homotopic maps) and $p_1 \sim - p_2$, hence $p_1 \sim - p_1$ which is impossible for an essential self-mapping on the two-sphere. More details of the given technique are found in \cite{Sjogren,\,Endo}.\hfill $\blacksquare$
\newline \newline

We have proven Proposition 1 on real matrix triples of order $q$. From this also follows, geometrically, that a complex matrix of order $n = 2m+1$ has a non-zero complex eigen-vector.

\head 2. De Medeiros' Vector Field Induced by an Endomorphism\endhead

We have just reviewed the Axis Theorem for $\CC^n$ in case the dimension $n = 2m+1$ is odd. The treatment we now give, of the full Axis Theorem for arbitrary $n > 0$ is perhaps easier as it does not require the topology of vector bundles, including stable Whitney sums, Grothendieck groups and the like. Instead, the methods used are from the result of H. Hopf (extending the same result for low dimensions due to Poincar\'e ), that for a closed manifold $M$, the ``index sum'' of any suitable smooth vector field on $M$ is pre-determined by $M$ (it is the integer Euler characteristic of $M$).

We try to avoid the business of computing Lefschetz numbers and Euler characteristics in favor of explicit constructions. The only variety that we actually deal with is the complex projective space $\CC P^n$, together with related spaces such as $\CC^{n+1}$.
This compact manifold $X = \CC P^n$ should be embedded smoothly into a (real) Euclidean space where it has a neighborhood $\Cal N$ that is a compact $2n$-manifold with boundary. Furthermore, $X$ will be a deformation retract of $\Cal N$.

In fact there exists a smooth (and can be chosen as holomorphic) vector field on $X$ that has exactly $n+1$ zeros, each occurring at a ``center'' of one of the canonical open sets that define the ``polar atlas'' of $X$. Each zero of this Hopf--Milnor field has index $ = 1$ by construction, so by the Poincar\'e--Hopf result noted above, any smooth vector field on $X$ must have a positive (non-zero) index sum.

On the other hand, an endomorphism $A$ of $\CC^{n+1}$ (linear self-mapping) is represented by an $n+1 \times n+1$ complex matrix $\hat{A}$. The rows and columns of $\hat{A}$ will be labeled $(0, 1, \dotsc, n)$ for the time being. The
matrix entries of $\hat{A}$ lead to a vector field on $\CC^{n+1}\, \backslash \,\{0\}$,
$$\Phi_{\hat{A}} = \sum_{j=0}^n\sum_{i=0}^n a_{ji} \,z_i \,\frac{\partial }{\partial z_j} \qquad  \mbox{(de Medeiros)}\quad .
$$
Here $\{z_j\}$ are the coordinate complex variables, though we sometimes also use $\{x_j\}$ as complex variables.  We note that this is exactly the kind of field, with ``linear'' coefficients at the tangent space basis $\left\{\frac{\partial}{\partial z_0}, \frac{\partial}{\partial z_1},\dotsc, \frac{\partial}{\partial z_n}\right\}$ that ``descends'' to a tangent vector field on $\CC P^n$. See Figure 2.

\bigskip
\centerline{\epsfbox{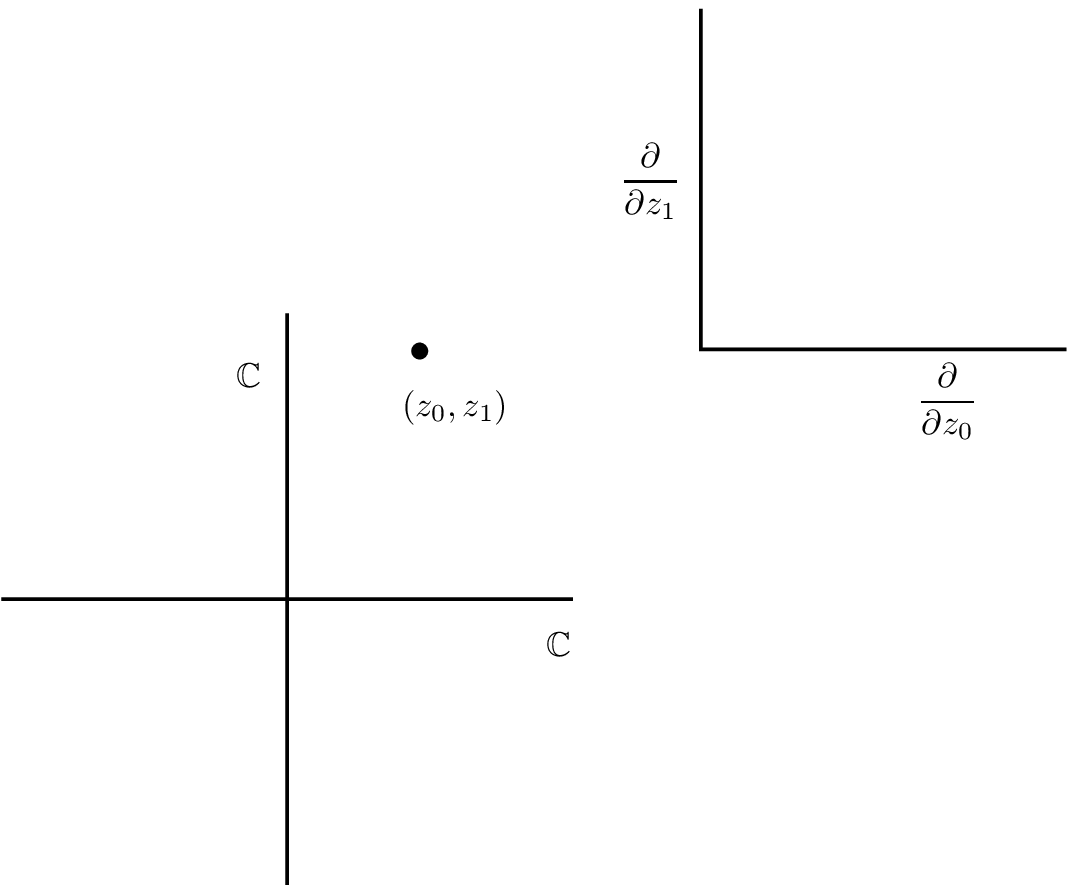}}

\medskip
\centerline{Figure 2. Complex Tangent Coordinates}

\bigskip
\noindent
{\bf Example  1}\quad Let $\hat{A} =
\left[\matrix
a & b & c \\
d & e & f \\
g & h & k
\endmatrix\right]
$ have complex entries, and $x, y, z$ be complex coordinates.

Now suppose that for certain values $x_0, y_0, z_0$ not all zero, there holds
$$\align
\phi_1 \equiv \quad & ax_0 + by_0 + cz_0 = \lambda x_0 \\
\phi_2 \equiv \quad & dx_0 + ey_0 + fz_0 = \lambda y_0 \\
\phi_3 \equiv \quad & gx_0 + hy_0 + kz_0 = \lambda z_0
\endalign
  $$
  where $\lambda \in \CC$.

  Then on $\CC^3$ we obtain
  $$\Phi_{\hat{A}} \left(x_0, y_0, z_0\right) = \lambda \left(x_0\frac{\partial}{\partial x} +
  y_0\frac{\partial}{\partial y} + z_0\frac{\partial}{\partial z}\right) \,\, .$$
  Thus for the chosen coordinate values, we have an outward- or inward-pointing tangent vector, called ``radial'', that descends or identifies to the zero (tangent) vector in $T_{[x_0:y_0:z_0]} \CC P^n$. See Figure 3.

  \bigskip
  \centerline{\epsfbox{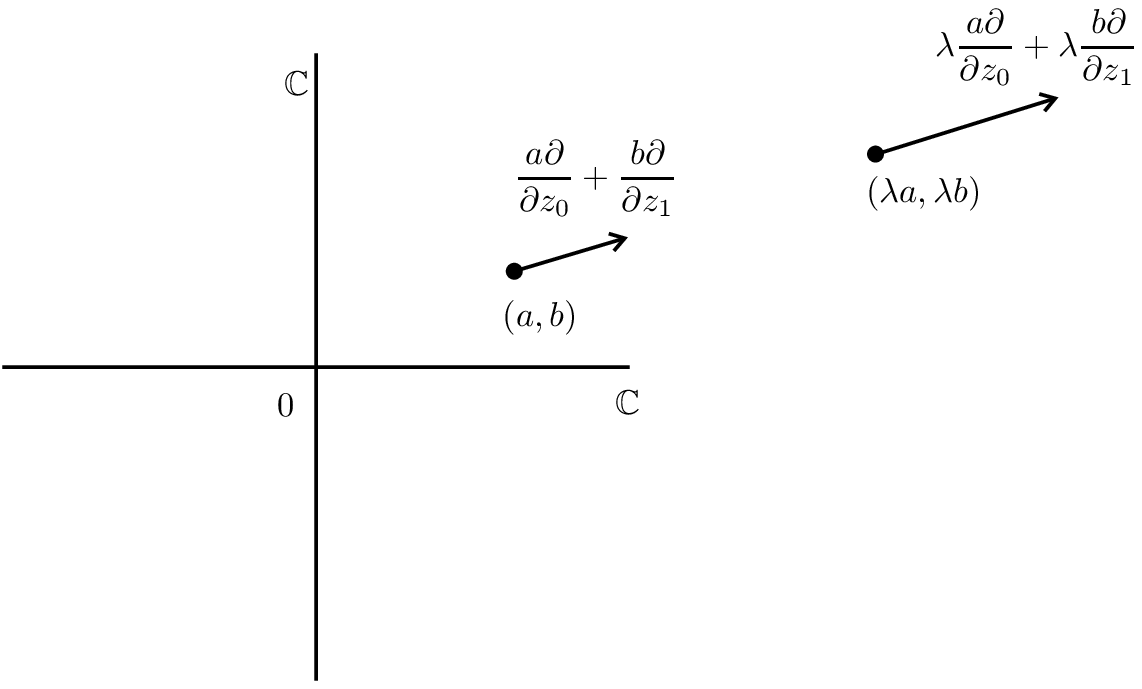}}

\medskip
\centerline{Figure 3. The Radial Vector Field on $\CC^2$}

\bigskip
  Conversely, if there is {\it no} solution $(x_0, y_0, z_0) \ne 0$ to
    $$\hat{A} \left(\matrix
    x_0 \\
    y_0\\
    z_0
    \endmatrix\right) = \lambda \left(\matrix
    x_0 \\
    y_0\\
    z_0
    \endmatrix\right)\quad ,$$
    we will {\it not} have
    $$\phi_1 (x_0, y_0, z_0) + \phi_2 (x_0, y_0, z_0) + \phi_3 (x_0, y_0, z_0)$$
    satisfying the radial property at any $(x_0, y_0, z_0) \in \CC^3$. In other words, the holomorphic vector field $\Phi_{\hat{A}}$ on $\Cal T \CC P^n$ (section of the tangent bundle) has no ``vanishing points''.
But the existence of such a vector field violates the Poincar\'e--Hopf theorem since we know that ${\Cal X}(X) = {\Cal X}(\CC P^n) > 0$ (the Euler characteristic).

    In the remainder of the paper, we fill out an argument that does not use numerical invariants from topology. Thus we intend to avoid computing or otherwise dealing with various Lefschetz numbers or Euler characteristics. It is not necessary to characterize the tangent bundle of projective space, though we do review the ``Euler sequence'' on background. The main tools used are the embedding of a compact manifold (into a real space), the tubular neighborhood, and Stokes' Theorem. Basic results about the calculus of differential forms are pointed to in the standard references. Particularly good is the text \cite{Edwards} on
many-variable Euclidean analysis.

\noindent
{\bf Theorem ``The Complex Axis''}\quad For any $n+1 \times n+1$ square matrix $\hat{A} = \frak M$ of complex scalars, there is an $n+1$-tuple $z_0, z_1,\dotsc, z_n \in \CC$, not all zero, and $\lambda \in \CC$, such that
$$\hat{A} \left[\matrix
z_0 \\
z_1 \\
\vdots \\
z_n
\endmatrix\right] = \lambda \left[\matrix
z_0 \\
z_1 \\
\vdots \\
z_n
\endmatrix\right]\quad .$$
Given any $\mu \in \CC$, $\mu \ne 0$, such a ``solution'' $\lambda$, $(z_0, \dotsc, z_n)$ leads to an equally valid solution $\lambda$, $(\mu z_0, \mu z_1, \dotsc, \mu z_n)$. The points $\mu \vec{z}$ constitute an Invariant Axis for $\hat{A}$, considered as missing the Origin. If the endomorphism $A$ is injective or surjective, $\hat{A}$ has an inverse and every such solution satisfies $\lambda \ne 0$.

\bigskip
\noindent
{\bf Example  2}\quad Consider ${\Cal C}$ of order 4 defined by
$${\Cal C} = \left[\matrix
0 & 0 & 0 & -\alpha \\
1 & 0 & 0 & -\beta \\
0 & 1 & 0 & -\gamma \\
0 & 0 & 1 & -\delta
\endmatrix\right],$$
so that
$$\lambda I_4 - {\Cal C} = \left[\matrix
\lambda & 0 & 0 & \alpha \\
-1 & \lambda & 0 & \beta \\
0 & -1 & \lambda & \gamma \\
0 & 0 & -1 & \lambda + \delta
\endmatrix\right]$$
and
$$\mbox{char}(\Cal C)= \mbox{det} \left(\lambda I - \Cal C\right) = \lambda^4 + \delta \lambda^3+\gamma \lambda^2 + \beta\lambda + \alpha \quad .$$
The matrix ${\Cal C}$ is companion to this characteristic (``secular'') polynomial, it has a solution of invariance (or Axis) if this polynomial equation has a root in the given ground field. For a contrasting philosophical point of view, see \cite{Axler}.

\head 3. The Explicit Tangent Bundle on $\CC P^n$\endhead

Given $p \in \CC P^n$, a vector $v_p \in T_p (\CC P^n)$ is determined by a ``curve'' $c(t)$ with $c(0) = p \in \CC P^n$, \, $0 \leq t < 1$. Consider now those vector fields on $\CC^{n + 1}-\{0\}$ which induce a field on $\CC P^n$ by projection to $\CC^{n+1}-\{0\}/\CC^*$.

The complex vectors must match up at $p \in \CC^{n+1}$ and at $\lambda p \in \CC^{n+1}$ for complex $\lambda \ne 0$.
But the curve that is equivalent to $c(t)$ at $\lambda p$ is just $\lambda c(t)$ which now has ``slope'' (or tangent vector) $\lambda c'(0)$ instead of $c'(0)$.

Holomorphic tangent vectors on $\CC^{n+1}$ can be written
$$v_p = \alpha_0 \left(x_0,\dotsc, x_n\right) \frac{\partial}{\partial x_0} + \cdots + \alpha_n \left(x_0,\dotsc, x_n\right) \frac{\partial}{\partial x_n}\quad ,$$
meaning
$$\alpha_0(p) \frac{\partial}{\partial x_0}+\cdots + \alpha_n(p) \frac{\partial}{\partial x_n}\quad .$$
The symbols $x_0,\dotsc, x_n$ in this instance denote complex variables. The {\it scaling} identity that we established as necessary for a well-defined descent to $\CC P^n$ means that $\alpha_j (\lambda p) = \lambda \alpha_j(p)$, so we may write
$$v_p = \beta_0 (p)\, \ell_0 \left(x_0,\dotsc, x_n\right) \frac{\partial}{\partial x_0}+\cdots +
\beta_n (p)\, \ell_n \left(x_0,\dotsc, x_n\right) \frac{\partial}{\partial x_n}\,\,\quad .$$
Here all $\ell_j \left(x_0,\dotsc, x_n\right) = \sum_{k=0}^n \gamma_{jk}\, x_k$ are linear forms with complex coefficients $\gamma_{jk}\in \CC$, and the $\{\beta_j\}$ satisfy a 0-homogeneous condition, namely that for all $p \in \CC^{n+1}$, there must hold $\beta_j (\lambda p) = \beta_j (p)$ when $\lambda \ne 0$. Thus the $\{\beta_j\}$ are all ``radially invariant''. When the base field is the field of complex numbers, this terminology should be interpreted with care.

We have seen that an endomorphism of $\CC^{n+1}$ induces such a ``de Medeiros vector field'' on $\CC^{n+1}$, which will project (descend) to a well-defined vector field on $\CC P^n$. Also the generating matrix $\hat{A}$ or $\frak M$ of the endormorphism gives rise through its {\it rows} $\left\{\left(\gamma_{jk}\right)\right\}$, $j = 0,\dotsc, n$, to the linear forms required for $v_p \in T_p \CC P^n$, varying holomorphically with $p$.

The characterization of the tangent bundle ${\Cal T}$ of $\CC P^n$ is often given in its relation to other vector bundles. In view of the equivalence of the concept ``vector bundle'', via the sheaf of sections, with the concept ``locally free sheaf'' on a variety, we may work with sheaf notation. See \cite{Gathmann} p\. 122-136 or \cite{Huybrechts}.

The important sheaves on $\CC P^n$ are $\Cal O$, $\Cal O (1)$ and $\Cal O (-1)$. A {\it sheaf} assigns to an open sub-set $U$ of the variety, a collection of functions into $\CC$, defined on $U$. The {\it structure sheaf} $\Cal O$ of $\CC P^n$ is first defined on {\it points} of $\CC P^n$, so that one may consider quotients $f/g$, where $f$ is a homogeneous polynomial in $(x_0,\dotsc, x_n)$ of some (total) degree $d$, and $g$ is another homogeneous polynomial of degree $d$, where we require the condition $g(p) \ne 0$. Then on an open set, we may define
$\Cal O (U) = \bigcap_{p\in U}{\Cal O}_p$.

The Serre twisting sheaf $\Cal O(1)(U)$ is given similarly to the above definition, where instead we look at quotients $f/g$ of homogeneous $f$ and $g$, with $g$ having total degree {\it one less} than does $f$. Analogously the {\it inverse Serre sheaf} $\Cal O(-1)(U)$ satisfies the provisions of $\Cal O (U)$ with the only change being that now the form $g$ should have degree {\it one greater} than that of $f$.

A foundational example is given on $\CC P^1$ with homogeneous coordinates $x_0, x_1$ that $\frac{1}{x_0} \in \Cal O (-1) (U)$, $U$ being a ``basic'' open set, namely $U = \left\{(x_0:x_1)\, \vert \,x_0\ne 0\right\}$. Here $f=1$ has degree $0$ and $g \equiv x_0$ has degree one.

The following ``Euler (short) exact sequence'' is connected to the well-known eponymous identity. We already saw the {\it radial} (or ``Euler'') {\it vector field}
$$\Delta_p = \sum_{j=1}^n x_j \,\frac{\partial}{\partial x_j}\,\,\quad.$$
This field may most easily be understood as defined on $\CC^{n+1}$. Also let $f$ be homogeneous of degree $d$. A typical ``campo vetorial de de Medeiros'' has degree $=1$. The classic {\bf Euler's identity} states that
$$\Delta_p f = d \cdot f(p)\, \quad.$$
The {\it Euler sequence} for the tangent bundle ${\Cal T}{(\CC P^n)}$ expresses
$$0 \to \Cal O \overset{i_*}\to\longrightarrow \oplus_{n+1} \Cal O(1) \overset\pi_*\to\longrightarrow {\Cal T} \to 0\,\,\quad.$$
Here the ``inclusion'' $i_*$ and ``projection'' $\pi_*$ are defined by

\roster
\item"(a)"
a ``function'' $h = f/g$ (0-homogeneous) maps to $h \cdot \Delta$. In coordinates for $\CC^{n+1}$, we pick $(x_0, \dotsc, x_n) \in \oplus_{n+1}\, \Cal O(1)$ which can then be scaled point-wise by the element $h$ of the structure sheaf $\Cal O$ belonging to $\CC P^n$. Thus $i_*\, h \simeq h \sum_{j=1}^n x_j\frac{\partial}{\partial x_j}$ as viewed from $\CC^{n+1}$.
\item"(b)"
having seen how any section of ${\Cal T}{(\CC^{n+1})}$ written $w = q \cdot \sum_{j=0}^n \ell_j \frac{\partial}{\partial x_j}$ projects to a derivation on $\CC P^n$, it now goes to $\pi_* w$ on ${\Cal T}{\CC P^n}$. Here $q$ is a 0-homogeneous function, and the $\{\ell_j\}$ are linear (1-homogeneous).
\endroster

Thus the kernel of $\pi_*$ is generated by the radial (Euler) field $\Delta$. This characterization of the tangent bundle of projective space, or more generally of a ``Grassmannian variety'', by means of its sheaves of sections, is handled in the usual texts, \cite{Harris}, \cite{Griffiths}, \cite{Hartshorne}, and \cite{Huybrechts} with some overlap.

For our purposes, a structure theorem such as the Euler sequence is included to provide context and is not really needed to prove the Axis Theorem. Given an endo-morphism, a `campo vetorial' of de Medeiros is constructed explicitly. We will see that the ``index'' of this vector field must equal that of another explicitly defined field on $\CC P^n$, the Milnor--Hopf field. This latter index for $n \geq 0$ is strictly positive.

If the original endo-morphism has {\it no} eigen-vector, its corresponding vector field has index equal to zero. This proves the Complex Axis Theorem, saying that such an endo-morphism and the resulting de Medeiros vector field cannot exist.

In the remainder of this report, we try to complete this demonstration without using too many globalistic tools besides methods of Calculus, mainly the Theorem of Stokes, interpreted as a result concerning the   co-bordism of hyper-surfaces.

\head 4. A Real Model of $\CC P^n$\endhead

We should be aware that our main topological space of interest, complex projective $n$--space, is a compact, complex manifold thus also retains the structure of a smooth (differentiable to all orders), compact real manifold (of dimension $2n$).

We will need to embed this manifold into some Euclidean space $\RR^N$ so that its given topology is the one induced as a subspace of this $\RR^N$. An embedding is a mapping that is proper, smooth, one-to-one, and is also injective via its differential, on the tangent space at any point.

Sharp estimates of a valid dimension for ambient embedding of a compact $k$-manifold $X$ were discovered by Whitney, and $2k+1$ is a safe number to take, so that we may embed $\CC P^n$ smoothly into $\RR^N$ for $N = 4n+1$.

Our purposes do not require a sharp estimate however. We note that Whitney's Theorem is often proved by selecting desirable projections from a dense set. We avoid such ``general position arguments'' which derive from measure theory, settling for a higher estimate on the dimension, that is easily derived, and constitutes the first step of the method of \cite{Whitney}.

The main facts about the complex projective space $X = \CC P^n$ can be explained through certain natural fibrations and canonical coverings. The equivalence relation $w \sim \mu w$, $\mu \in \RR^+$ gives the projection (a) below, to the compact real variety $\SS^{2n+1}$, which may be considered the locus of $|z_0|^2+\cdots +|z_n|^2 = 1$. Recall that an equatorial circle $\SS^1$ is a Lie group parametrized by angles $0 \leq \theta < 2 \pi$. Any complex line $L \subset \CC^{n+1}$ intersects $\SS^{2n+1}$ in such a great circle, so if the points of the circle are identified, we have the
coordinates of a line: in other words, an element $[L] \in \CC P^n$. This yields the projection (b) in the expression
$$
\CC^{n+1}\backslash\{0\}\overset{(a)}\to\longrightarrow \SS^{2n+1} \overset{(b)}\to\longrightarrow\CC P^n \quad .\tag{$\text{D}$}
$$
Since $\SS^{2n+1} \subset \RR^{2n+2}$ is closed and bounded in a metric space, it is compact, and (b) is a continuous identification mapping\footnote{Because (b) is a ``projection''.}, so $\CC P^n$ is compact. Furthermore, $\SS^{2n+1}$ is Hausdorff with the endowed topology, and the ``action'' of $\SS^1$ on $\SS^{2n+1}$ is continuous, i.e. $\alpha (\theta, \vec{y}) = e^{i\theta}\vec{y}$ where $\alpha: \SS^1 \times \SS^{2n+1} \to \SS^{2n+1}$.

The great circles on an odd sphere are disjoint. In any case we have that for $\SS^{2n+1} \overset{(b)}\to\longrightarrow \SS^{2n+1} / \,\SS^1 = \CC P^n$, the mapping (b) is continuous and open. As a compact Hausdorff space, it is also {\it normal} (disjoint closed sets may be separated by disjoint open sets that contain them). For the classical ``Hopf bundle'' when $n = 1$, a rough depiction is given in \cite{Hatcher} p\. 377.

\bigskip
For the purpose of forming a suitable embedding of $\CC P^n$ into Euclidean space, we review an atlas for its structure as a complex manifold. For $j = 0, \dotsc, n$ let
$$U_j = \left\{\left(z_0: \cdots : z_j: \cdots : z_n\right)\left| \,z_j \ne 0\right. \right\}\quad .
$$
Here as usual, $p = (z_0:\cdots : z_n)$ means the projection $\pi: \CC^{n+1} \to \CC P^n$ applied to $(z_0,\dotsc, z_n) \in \CC^{n+1}$. Then $\varphi_j: U_j \to \CC^n$ defined by $(z_0: \cdots : z_n) \mapsto \left(\frac{z_0}{z_j},\dotsc, \frac{z_n}{z_j}\right)$ is continuous as follows from the fact that $\pi: \CC^{n+1} \to \CC P^n$ is continuous. Also $\varphi_j$ is {\it open} as follows from the fact that $\pi$ is open. Summarizing, $\varphi_j$ is also bijective, hence a homeomorphism.

To complete the complex atlas, for $j < i$ say, we calculate transition mappings $\psi_{ij} = \varphi_j \circ \varphi_i^{-1}$, $\varphi_i (U_i \cap U_j) \to \CC^n$ as
$$\psi_{ij} (y_1,\dotsc, y_n) = \left(\frac{y_1}{y_j},\dotsc, \frac{y_{j-1}}{y_j},\frac{y_{j+1}}{y_j},\dotsc, \frac{y_i}{y_j}, \frac{1}{y_j},\frac{y_{i+1}}{y_j},\dotsc, \frac{y_n}{y_j}\right)\quad .$$

These transition ``functions'' are biholomorphic, so viewed as $\psi_{ij}: V_{ij} \to \RR^{2n}_w$, they are diffeomorphisms to the image. The open set $V_{i j} \subset \RR^{2n}_z$ arises from $U_i \cap U_j$ by realization.

To find an embedding into some $\RR^N$ of a smooth compact manifold $X$, such as the realization of $\CC P^n$, we work in the category of (infinitely) smooth (real, abstract) manifolds. The standard construction uses ``bump functions'' such as $B: \RR \to \RR$, smooth and satisfying $B(x) > 0$ for $|x| < 1$ and $B(x) = 0$ for $|x| \ge 1$. The typical exemplar which is used is $B(x) = e^{-(x^2+1)/(x^2-1)^2}$ for $|x| < 1$. This kind of bump should not exist in a holomorphic or real analytic category. Moving to a higher dimension, there are bump functions Taylor-made for manifolds.

\bigskip
\noindent
{\bf Lemma}\quad Let there be given a smooth abstract manifold $M$ with $K$ compact contained in an open $U \subset M$. Then there exists a smooth function $g: M \to \RR^+$ satisfying $g(x) > 0$ for $x$ in $K$, and $g(y) = 0$ for $y \in M \,\backslash \,U$.
Thus the ``support'' of $g$ lies within $U$, and $g$ {\it bumps} to a positive value on the compact $K$.

\noindent
{\it Proof}\quad See \cite{Bredon}, section II.10.\hfill $\blacksquare$

\bigskip
Consider again the canonical atlas $\{(\varphi_j, U_j)\}$ of $X = \CC P^n$. Since $X$ is normal as a topological space, the ``shrinking lemma'' applies and we have closed $K_j \subset U_j$, $j=0,\dotsc, n$ so that $\{\mbox{Int}(K_j)\}$ {\it also} covers $X$. From the above Lemma, one may show that there exist $\{\lambda_j: X \to \RR\}$, smooth so that $\lambda_j (x) = 1$ for $x \in K_j$, with $\mbox{supp}(\lambda_j) \subset U_j$. Next, we define smooth functions with values in $\RR^{2n}$
$$\sigma_j(p)  =
\left\{\matrix
\lambda_j(p)\,\varphi_j(p), & \mbox{for $p\in U_j$ \quad , }\\
\vec 0,  & \mbox{for $p\notin U_j$\qquad \quad.}
\endmatrix \right.
$$
We recall that the number of charts $U_j$ in the atlas is $n+1$. Putting the constructions together allows the definition of $\gamma: X \to \RR^{2n(n+1)} \times \RR^{n+1}$ by means of $\gamma (p) = \left(\sigma_0 (p), \dotsc, \sigma_n (p), \lambda_0 (p),\dotsc, \lambda_n (p)\right)$.

From the fact that the $\{\varphi_j\}$ are diffeomorphisms, their induced tangent mappings are injective, so also is $\gamma_{* p}: TX_p \to T_{\sigma(p)} \left(\RR^{2n(n+1)}\right) \times T_{\lambda(p)}\left(\RR^{n+1}\right)$, hence $\gamma$ is confirmed to act as an immersion mapping.

Since $X = \CC  P^n$ is compact as seen by the projection (b), the mapping $\gamma$ is automatically proper.
Similarly, $X$ is orientable. This may also be deduced from the decomposition of $X$ into real cells of even dimension. The CW ``attaching'' procedure is detailed in the Introduction, with graphic formulas. The attaching process will map the now cell at its spherical boundary to the whole space $\CC P^{k-1}$.
The full fiber space maps, according to the Hopf fibration, to the ``highest'' space $\CC P^{k-1}$ that has recursively been constructed up to this point.

Finally we would like to demonstrate that $\gamma$ is a one-to-one mapping. Supposing $\lambda_j (p) = \lambda_j(q)$. But for some $i \in [0, n]$ we observe that $p \in \mbox{Int}(K_i)$ where follows $\lambda_i (p) = 1$.
Then $\sigma_i (p) = \lambda_i(p) \varphi_i(p) = \lambda_i (q) \varphi_i (q) = \sigma_i (q)$. But $\varphi_i$ is a homeomorphism from $U_i$ to its image in $X$, hence one-to-one, therefore $p=q$ as expected.
Using the count of atlas members, we have arrived at an explicit mapping $\gamma : X \to \RR^{(n+1)(2n+1)}$.

Besides the well-known reduction of the ambient dimension to $4n+1$ or even $4n$, due to \cite{Whitney}, there are several ways to reduce the dimension for our specific $2n$-manifold $X$, without employing general position arguments. Instead of embedding an arbitrary compact manifold, we concentrate on the projective spaces of interest.
For example, \cite{James} used an algebraic construction, taking into account the tangent space, that gives embedding dimensions for real, complex, and quaternionic projective spaces. This yields an embedding of $\CC P^n$ into $\RR^{4n-1}$, stronger than Whitney's general result, without the use of Sard's Theorem.

R\. J\. Milgram has exposited a way to treat immersions of projective spaces by means of an explicit collection of real skew-symmetric matrices. From this technique there arise embeddings that often exhibit a sharper embedding dimension for $\CC P^n$. See \cite{Mukherjee}, \cite{Steer}, \cite{Milgram}.

Of course, it is of interest to consider {\it isometric} embeddings. Here we mean that $\CC P^n$ be endowed with the Riemannian metric obtained by means of our Hopf projection $\SS^{2n+1}\overset{(b)}\to\longrightarrow \CC P^n$, starting with the standard ``round'' metric on $\SS^{2n+1} \subset \RR^{2n+2}$. Straightforward application of Nash's Embedding Theorem yields an ambient dimension (for some isometry) of $N = 6n^2+11 n$. Based on important work of \cite{G\"unther}, one may obtain an improved ambient dimension
$$N(n) = \left\{\matrix
45, & n < 5, \\
2(n^2+3n),  & n \geq 5
\endmatrix\right.$$
for $\CC P^n$ isometrically embedded into $\RR^N$. That is, even under the ``isometry'' condition, we obtain an embedding dimension that is computable and not inordinately large. Whichever embedding is chosen, continue to write it as $\gamma: X = \CC P^n \to \RR^N$. See \cite{Konnov} and \cite{Lu, Isometric} for further background on isometric embedding.

\head 5. Extension of a Non-degenerate Field to the Tubular Neighborhood of Embedded $X$\endhead

We wished to prove that a ``campo vetorial'' (field) $w$ according to de Medeiros, living on $\CC P^n$, $n \geq 0$, must have a vanishing point (similar to as in the ``Hedgehog Theorem'', see Milnor's article \cite{CUE}). Since the Euler characteristic of such a manifold $X$ is non-zero, our result follows quickly from the Poincar\'e--Hopf Theorem, see \cite{G-P}. 

\bigskip
As pointed out in the Introduction, the short proof of the Complex Axis Theorem by de Medeiros is essentially a homotopy argument, and for $\mbox{id}:\CC P^n \to \CC P^n$, exploits the formula
${\Cal L}(\mbox{id}) = {\chi} (\CC P^n) = n+1$.

For the reasons listed in the Introduction, we propose an alternative, which is actually a co-bordism argument that does not require any explicit numerical invariant to be calculated. Co-bordism, actualized through Stokes' Theorem, allows the comparison of any two non-degenerate fields via their indices. Again, the ``index'' measures the required quantity of vanishing points of the field.

We invoke a relationship between portions of the boundary of an $N$-manifold. But there exist fields arising from an $(n+1)$-endomorphism, in other words ``campos vectoriais de de Medeiros'' with each zero occurring in a positive sense. This shows that such a vector field is impossible, and from there using the reciprocal relationship between ``secular polynomial'' and ``companion matrix'', we arrive at another proof to the Fundamental Theorem of Algebra.

The guiding light in this program is the Lemma of Hopf, which was a key step in his generalization of Poincar\'e's Theorem on vector fields that possess isolated zeros. For our purposes we may pare down the proofs and definitions required, by restricting ourselves to vector fields on $X$ whose zeros are {\it non-degenerate}. Thus if $v(z_0) = \vec{B} \in \RR^{2n}$, we are assuming that the differential $dv$ at $z_0$ is non-singular (as an endomorphism of $\RR^{2n}$, or a matrix). Equivalently there is an open ball $B_{z_0} \subset X$, small in the Fubini--Study metric, which maps diffeomorphically under $v$ to its image in $\RR^{2n}$. These stipulations simplify matters by constraining the Poincar\'e--Hopf {\it index} of $v$ at $z_0$ either to be $+1$ or $- 1$ (orientation-preserving diffeomorphism or not).

\bigskip
\noindent
{\bf Hopf\rq s Lemma}\quad Suppose that $M \subset \RR^N$ is a compact $N$-dimensional smoothly embedded manifold-with-boundary. Thus $\partial M$ is a smooth $(N-1)$-dimensional sub-manifold of $\RR^N$ with one or more component. Let $V$ be a smooth vector field having (finitely many) non-degenerate zeros only, such that ``$v$ is outward-pointing on $\partial M$''. We then have a Gau{\ss} (normal) map $G: \partial M \to \SS^{N-1}$. Summing the indices at the zeros of $v$, we obtain
$$\mbox{degree $G$} = \sum_{\text{$z \in$  zeros $(v)$}} \iota_z\quad .$$
\bigskip
\noindent
{\bf Corollary 1}\quad Any two such vector fields on $M^N$ have the same index sum.
\bigskip
We shall prove Hopf's Lemma after giving some treatment of the Tubular Neighborhood concept. For a fuller exposition see \cite{Spivak}, Vol\. I, Chap\. 9 (Addendum). The ``persistent rumor'', given in square brackets on p\. 299 of \cite{Spivak}, Vol\. III, that when $m$ is even, the degree of the normal map for a hyper-surface $L \subset \RR^m$ equals one-half $\chi(L)$ should be considered in the light of remarks by Milnor, \cite{TFDV}, p\.~86. In particular, for the standard embedded $\SS^1 \subset \RR^2$, one sees that the Gau{\ss} map ought to have (Brouwer) degree equal to $1$.

Working with a (de Medeiros) vector field $v$ on $X^{2n} = \CC P^n$ using Hopf's Lemma, we need to ``fatten'' $X$ to a sub-manifold $W^N\subset \RR^N$, the ambient Euclidean space, and extend $v$ appropriately.
If the compact manifold $X$ of dimension $q = 2n$ is smoothly embedded into $\RR^N$, the tangent ``plane'' of dimension $q$ is defined at each $x \in X$, denoted $T_x$. We may construct a space $\frak N(X)$ defined by
$$\left\{(x, y)\,|\, x \in X,\, y \in \RR^N\;\text{with}\,\;y \cdot t = 0\;\text{for all}\,\;t \in T_x\right\}\quad .$$
This is the ``normal bundle'' to $X \subset \RR^N$, and we have
$\pi : \frak N(X) \to X$ given by projection to the first coordinate. Furthermore, $x + y$ is always defined in $\RR^N$, so there is a mapping $\theta: \frak N(X) \to \RR^N$ by $\theta (x, y) = x+y\,$.

Similarly define $\frak N_{\epsilon}(x) $ by requiring $\|y\| < \epsilon$, and put
$$W_{\epsilon} = \left\{y \in \RR^N: \text{inf}\left(\|y - x\| < \epsilon\quad\text{over}\quad x\in X\right)\right\}\quad .$$

\bigskip
\noindent
{\bf Tubular Neighborhood Theorem}\quad Under the stated hypotheses, there exists $\epsilon > 0$ such that $W_{\epsilon} = \theta (\frak N_{\epsilon}(x))$ comes about from a diffeomorphism.

\bigskip
\noindent
{\bf Remark}\quad By the normal bundle con\-struc\-tion, we see that $W_{\epsilon}$ is a smooth manifold-with-boundary of dimension $N$. Furthermore $\pi \cdot \theta^{-1}: W_{\epsilon} \to X$ is a deformation retraction, so that $\pi \circ \theta^{-1} \circ \gamma: X \to X$ is defined, given the embedding $\gamma: X \to \RR^N$. We would also know that $\gamma \circ \pi \circ \theta^{-1}: W_{\epsilon}\to W_{\epsilon}$ is smoothly homotopic to the identity of $W_{\epsilon}$ through a family of diffeomorphisms.

\bigskip
\noindent
{\it Proof of Tubular Neighborhood Theorem}\quad (For details see \cite{Bredon}, Chap\. II\.) One observes that the tangent bundle $\Cal T(\RR^N)$ splits by direct sum into the parts tangent to and normal to $X$ in $\RR^N$. Examining the mapping $\theta$ on normal and tangent directions, we deduce that for the differential $d \theta :{\Cal T} (\frak N(X)) \to {\Cal T}(\RR^N)$ defined by $d \theta_x: T_{x, 0}\, \frak N \to T_x (\RR^N) \simeq \RR^N$, $\,x \in X$, has full rank. Thus for $(x, 0) \in \frak N(X)$ we have by the Inverse Function Theorem, that $\theta$ is a diffeomorphism on a neighborhood of $X$, $V_X \subset \frak N(X)$. Now the compactness and sequential compactness of $X$ show that $\epsilon > 0$ can be found so that $\theta : \frak N_{\epsilon}(X) \to \RR^N$ is a diffeomorphism on some neighborhood of each point, and also $\theta$ is one-to-one. Finally we made the assertion that the image under $\theta$ of $\frak N_{\epsilon}(X)$ is just the metrical $\epsilon$-neighborhood of $X$ lying within $\RR^N$. But this follows from the fact that for small $\epsilon > 0$, given $\|x - y\| < \epsilon$ for all $x \in X$, the $\hat{x}$ value where this distance attains a minimum yields: $\hat{x}-y$ is normal to $X$ (it is orthogonal to any vector of $T_{\hat{x}}(X)$). \hfill $\blacksquare$

\bigskip
To apply Hopf's Lemma, we need, given a non-degenerate vector field $v$ on $X$, to construct a compatible field $w$ defined on $W_{\epsilon}$. We may furthermore regard $W_{\epsilon}$ as invested with its natural boundary in $\RR^N$, so that it is seen to be a manifold-with-boundary of dimension $N$, and the expression ``$w$ is outward-pointing'' makes sense. In fact the new field $w$ will have zeros {\it only} lying on $X$, the {\it same} zeros with the {\it same} indices as has $v$.

We may write the projection $\pi: W_{\epsilon} \to X$ with mild abuse of notation. Given $q \in W_{\epsilon}$, the extended vector field should be defined as $w(q) = [q - \pi(q)] + v(\pi(q))$. The terms in square brackets lie in $\RR^N$ as a result of the embedding $W_{\epsilon} \subset \RR^N$. The final term lies in $\RR^N$ through definition of a vector field on $X \subset \RR^N$. The vector $w(q)$ of the new field should be considered as ``based at $q \in W_{\epsilon}$''. We already saw why $q - \pi(q)$ is orthogonal to $T_{\pi(q)}(X)$.

A calculation shows that $w(q)$ is outward-pointing on $\partial W_{\epsilon}$. Indeed, if $\rho(q) = \|q - \pi(q)\|$, the gradient of $\rho$ at $q$, \,$\text{grad}\,\rho = \displaystyle{\sum_{i=1}^N \frac{\partial \rho(q)}{\partial x_i} \cdot \frac{\partial}{\partial x_i}}\,$ equals $\,2(q - \pi(q))$. Thus, the unit normal to the surface $\partial W_{\epsilon} = \rho^{-1}(\epsilon)$ becomes $h(q) = (\text{grad}\,\rho)/\|\text{grad}\,\rho\| = \frac{1}{\epsilon}\left\{q -\pi(q\right\}$, and the dot product $w(q)\cdot h(q) = \epsilon > 0$ on $\partial W_{\epsilon}$. Next suppose that $w(q) = 0$ for some $q \in W_{\epsilon}$ (which must be in the interior). If $q - \pi(q)$ is non-zero, it is also orthogonal to $v(\pi(q))$, so $w(q) \ne \vec 0$. Hence $v(\pi(q)) = v(q) = \vec 0$. All zeros of $w$ in $W_{\epsilon} \simeq \frak N_{\epsilon}$ lie in $X$ and are also the given zeros of $v$. See Figure 4.

\bigskip
\centerline{\epsfbox{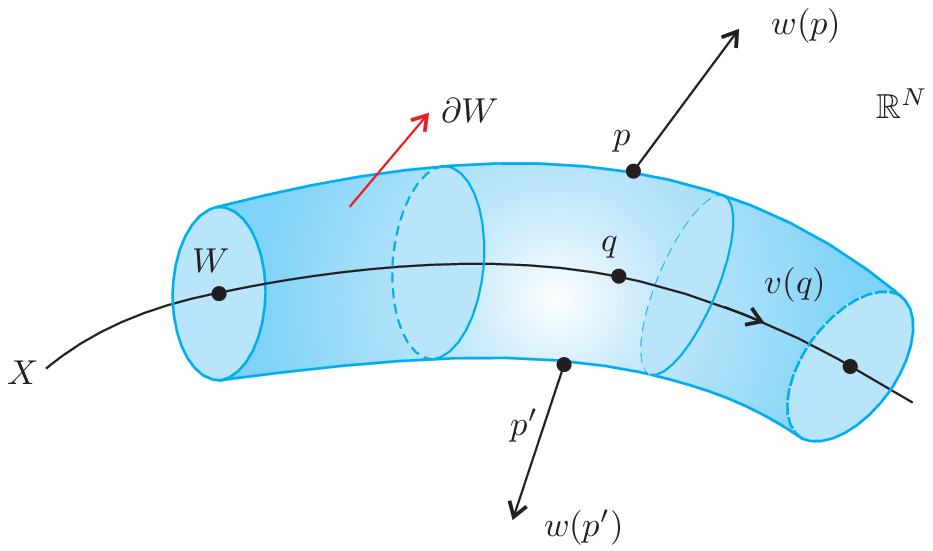}}

\medskip
\centerline{Figure 4. Extension of Tangent Vector Field on an}
\smallskip
\centerline{Embedded Manifold to a Tubular Neighbourhood}

\bigskip
Thus we have constructed a compatible vector field $w$ on $W_{\epsilon}$, extending the original field $v$ defined on $X$. Finally, to compute the index (``mapping degree'') of the vector field $w$ at a point $q \in X = \pi (W_{\epsilon})$, it is sufficient to consider a vector $\kappa_{\top} \in T_q(X)$ and a vector $\kappa_{\bot}\in \frak N(q)$, with $\kappa_{\top} \cdot \kappa_{\bot} = 0$ holding. The collection of these vectors $\{\kappa_{\top}, \kappa_{\bot}\}$ span all of $T_q(\RR^N) \simeq \RR^N$. We have the following action of a differential form on a vector, thus
$$\align
dw(q)\, \kappa_{\top}\, & = \,dv(q)\, \kappa_{\top} \\
dw(q)\, \kappa_{\bot}\, & = \, \kappa_{\bot}\qquad \quad\text{(the identity linear mapping)\quad .}
\endalign
  $$
  From the determinant of a ``block diagonal'' matrix we obtain
  $$\text{det}\,dw (q) = \text{det}\,dv (q)$$
or, that $dw (q)$ is non-singular and orientation-preserving just when $dv(q)$ also has these properties. In other words, $\text{index}\,w(q) = \text{index}\,v(q)\,$, where $q$ is a ``zero''. This could also be shown without the use of determinants. These index values in the case of non-degenerate $v$ certainly are equal to $\pm 1$.

  \head 6. The Volume Form and Hopf's Lemma\endhead

  The work we have completed in Section 5\. facilitates certain deductions from Hopf's Lemma.

\noindent
{\bf Corollary 2}\quad  Given two vector fields $v_a$, $v_b$ on $X^{2n} = \CC P^n$  with all zeros non-degenerate, their index sums $\sum \iota \,v_a,\, \sum \iota \,v_b$ are the same, equal to the degree of the Gau{\ss} map $G: \partial W_{\epsilon} \to \SS^{N - 1}$ for a suitable common tubular neighborhood $W_{\epsilon}$ of $X = \CC P^{n}$.

In particular, if $v_a$ has no zeros at all, and at some zero $p_1 \in X$ of $v_b$, the index is positive $(= +1)$, then at some other zero $p_2$, $p_1 \ne p_2$, $\iota_{v_b}(p_2)$ must be negative $(= - 1)$.

We may limit our interest in the Volume Form to the ambient space $\RR^N$ and to a closed hyper-surface $M \subset \RR^N$ (of dimension $N - 1$, but not necessarily connected).
An $N$-form $\eta$ may be integrated over any oriented Jordan region ${\Cal J}$ of $\RR^N$. Or we may denote
$I({\Cal J}) = \int_{\Cal J} f \cdot \eta$ to be the {\it integral} of $f: \RR^N \to \RR$, $f$ being continuous on disjoint Jordan subsets of ${\Cal J}$ that exhaust ${\Cal J}$. In case for a constant mapping $f(x) = 1$, $x \in {\Cal J}$, and $I({\Cal J})$ always adds up to the (oriented) Jordan content of ${\Cal J}$, we say that $\eta$ is a {\it volume form} on $\RR^N$. It is well-known that $\eta$ has to equal $dx_1 \wedge \cdots \wedge dx_N$, see \cite{Flanders}.

\bigskip
\centerline{\epsfbox{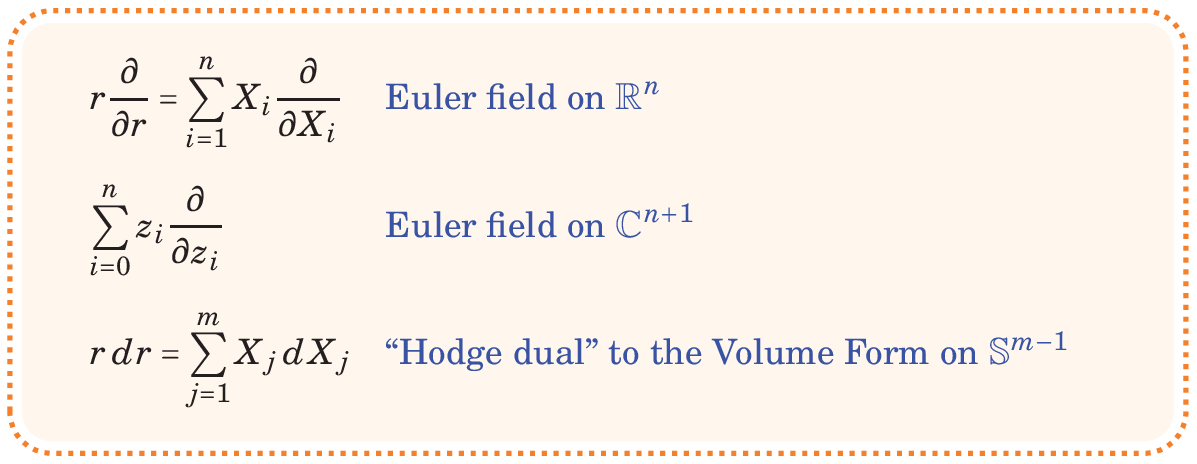}}

\medskip
\centerline{Figure 5. Rotationally Invariant Fields and Forms}

\bigskip
The other volume form that we consider explicitly lives on $\SS^{N - 1}$, the unit sphere which inherits its metrical properties through the canonical Euler embedding. See Figure 5 for relations between ``radial quantities''.
To construct the canonical volume form on $\SS^{N - 1}$, it is helpful to review the Hodge duality of Euclidean forms. On $\RR^N$ define
$$\align
*\, & dx_1 = dx_2 \wedge \cdots \wedge dx_N \\
*\, & dx_2 = (- 1) \cdot dx_1 \wedge dx_3 \wedge \cdots \wedge dx_N \\
&\vdots \\
*\, & dx_N = (- 1)^{N+1} \cdot dx_1 \wedge dx_2 \wedge \cdots \wedge dx_{N - 1} \quad . \\
\endalign  $$

The sign arises according to whether the indices of $\left(x_j,\, x_1, \,x_2,\,\dotsc,\, \hat{x}_j,\, \dotsc,\, x_N\right)$ form an {\it even} or {\it odd} permutation of $[N]$. See the report of \cite{Dray} for detailed calculations.

To construct the preferred Volume Form on $\SS^{N - 1}$ we start with the differential $1$-form $\tau = \displaystyle{\sum_{i=1}^N}\, x_i \,dx_i$ on $\RR^N$. This appears to be {\it dual} to the Eulerian radial vector $\displaystyle{\sum_{i=1}^N x_i \,\frac{\partial}{\partial x_i}}$ which we have already seen in a complex guise. In any case, a calculation shows that $\tau = r\, dr$ in spherical coordinates and hence this form is invariant under the group of rotations $\bold S \bold O (N)$. Therefore the Hodge ``star duality'' $(N - 1)$-form $*\,\tau$ must also be rotation invariant. In other words, a Jordan patch ${\Cal J}\subset \SS^{N-1}$, rotated to $\gamma({\Cal J}) \subset \SS^{N - 1}$ will have the same $(N - 1)$-area as before. By the rostered formulas above, we obtain
$$\omega = *\,\tau = \sum_{i=1}^N (- 1)^{i - 1} x_i\, dx_1 \wedge \cdots \wedge \hat{dx}_i \wedge \cdots dx_N \qquad .$$
One may feel confident to have found the canonical $(N - 1)$-form on $\SS^{N - 1}$ up to a constant factor. But by Stokes' Theorem,
$$\int_{\SS^{N-1}} \omega = \int_{\BB^N}d\,\omega = \int_{\BB^N} N \,dx_1\wedge \cdots \wedge dx_N = N \cdot {\text Vol}(\BB^N)\quad .$$
Here $\partial \BB^N = \SS^{N - 1}$. But this quantity $N \cdot {\text Vol}(\BB^N)$ is the correct value for the ``area'' of the unit $(N - 1)$-sphere, as shown in popular expositions such as \cite{Folland} and \cite{XWang}. Therefore the scaling of $*\,\tau$ we gave must be correct. Alternatively, in spherical coordinates with $\theta_1 \in [0, \pi], \,\dotsc,\, \theta_{N - 2} \in [0,\pi],\, \theta_{N - 1} \in [0, 2\pi]$, we would write
$$\omega = d\SS^{N - 1} = \sin^{N - 2} (\theta_1)\, \sin^{N - 3}(\theta_2) \cdots \sin (\theta_{N - 2}) \cdot d\theta_1\wedge \cdots \wedge \theta_{N - 1}\quad ,$$
as derived in \cite{Blumenson}.

We return to a consideration of ``Hopf's Lemma''. In fact this result is valid more generally than we have stated it. The vector field $v$ could possess {\it degenerate} zeros (on the $N$-manifold, with boundary, $M$) as long as these are finite in number. Thus the vector field near a zero $x_0 \in M$ would not need to satisfy $\iota_v \,(x_0) = \pm 1$. It would not necessarily give a diffeomorphism on a neighborhood of $x_0$.

The more specialized ``non-degenerate'' case is however all that we need to complete the Complex Axis Theorem, and to confirm the utility of the mapping $\Lambda: \text{End}(\CC^{n+1}) \to \text{Vect}(\CC P^n)$, due to von Sohsten de Medeiros. The {\bf Complex Axis Theorem} states, given $\frak M$ a complex $n+1 \times n+1$ matrix, that $\Lambda (\frak M)$ is a vector field that must have a vanishing point (``zero'') on $\CC P^n$.

We will see that this ``vector field'' approach amounts to a bordism argument. The zeros of (non-degenerate) $v$ give rise to neighborhoods of $M$ and their (spherical) boundaries. The oriented sum of those spherical boundary components is co-bordant to the remaining boundary component $\partial M$, so that their indices (mapping degrees) add up to the degree of the Gau{\ss} (normal) mapping $G: \partial M = \partial W_{\epsilon} \to \SS^{N - 1}$.
In our application of Hopf's Lemma, in fact we need only the case of ``positively oriented'' zeros. Thus, an index $\iota (x_0)$ for a zero $x_0$ of $v$ can be taken  as equal to $+1$.

The striking proof of Complex Axis by \cite{de Medeiros} involves a homotopy argument. Generally, constructions via homotopy are more delicate than the technique of looking at co-bordant hypersurfaces. It is of interest to write down proofs of important theorems using the most primitive implements possible. This is one reason to propose a ``bordism'' modification to the method of A\. v\. S\. de Medeiros.

\bigskip
\noindent
{\it Proof of Hopf's Lemma (completed)}\quad Consider $v$ on $M^N$ as given, with zeros in quantity $k$, $q_1 \in U_1$, $q_2 \in U_2,\,\dotsc,\, q_k \in U_k$ contained in small balls ($N$-disks) as indicated. 

\bigskip
Let $Q = M \, \backslash \,\, (\,U_1 \cup \,\cdots\, \cup U_k\,)$, a manifold-with-boundary of dimension $N$, where the smooth vector field $v$ is supported and is non-vanishing. Thus the normalized mapping $G(x) = \frac{v}{\|v\|}$ is well-defined and takes $Q$ to $\SS^{N - 1}$. We have constructed some canonical volume forms and now we use the $(N - 1)$-differential form $\omega$  toward our bordism situation, relating $\partial M$ to  $\,\partial U_1 \,\cup\, \cdots \,\cup\, \partial U_k$ by means of Stokes' Theorem.

First of all, the axiomatics of ``forms'' yield a computation for the exterior derivative, namely $dG^*(\omega) = G^*(d\omega)$ from \cite{Flanders}, p\. 24.  Here $G^*$ is the mapping induced by $G$ on differential forms. By Stokes,
$$\int_{\partial Q} G^*(\omega) = \int_Q d \,G^*(\omega) = \int_Q G^* d \omega\quad .$$
We indicate why the final integral must equal $0$. Examine $$\int_c G^* d\omega \quad ,$$ integration over an $N$-simplex $c$, geometrically contained as $c \subset Q \subset \RR^N$. As explained in \cite{Flanders} p\. 73, one may take the (co-variant) mapping $G_*$ on chains, getting an equal quantity $\int_{G_* c} d\omega$. Since $G$ maps $Q$ to the lower-dimensional space $\SS^{N - 1}$, the chain $G_* c$ must be degenerate (it is geometrically of dimension $N - 1$). But the canonical volume form $\omega$ yields $d\omega$ as an $N$-form (a multiple of $dx_1 \wedge \cdots \wedge dx_N$, which is naturally defined on $\RR^N \,\backslash \,\{0\}$). Hence the right-most integral expression above must equal zero. 
\bigskip
Now degree$\,(G)$ on $\partial M$ is just the integer coefficient $\mu$ in the defining formula
$$\int_{\partial M} G^*(\omega) = \mu \cdot \int_{\SS^{N - 1}} \omega = \mu \cdot {\text Area} (\SS^{N - 1})\quad .$$
Letting $N = 2a$ or $2a+1$ for a natural number $a$, the $(N - 1)$-dimensional ``area'' $A(\SS^{N - 1})$ is given by
$$A(\SS^{N-1}) = \frac{2^{a+1} \pi^a}{1\cdot 3\cdot 5 \cdots (N - 2)}\quad \text{for $N$ odd}\,$$
$$A(\SS^{N - 1}) = \frac{2 \cdot \pi^a}{(a - 1)!}\quad \text{for $N$ even}\qquad ,$$
when $N \geq 1$, so these hyper-areas have values greater than $0$. When
$N - 1 = 0$ we have a counting measure on the sphere, so $A (\SS^0) = 2$. For various ways to calculate $A( \SS^{N - 1})$ and $V(\BB^N)$, area and volume of sphere and ball respectively, see \cite{XWang} and references therein.

Getting back to the integral of the normal map, similar integrals over the components $\partial U_1,\dotsc , \partial U_k$ are {\it additive}, leading to
$$\frac{1}{A(\SS^{N-1})} \int_{\partial Q}G^*(\omega) = \text{degree}(v\;\text{on}\;\partial M) -
\sum_{i=1}^k \iota_v (q_i) = 0\quad .$$
Thus in particular, the sum of indices of the (non-degenerate) field chosen on $M$ does not depend upon which such field is chosen. All smooth, non-degenerate fields yield the same index sum, an integer that can be computed from $M$, see \cite{Flanders} Section 6\.2. As we are emphasizing the case of $M$ as the tubular neighborhood $W_\epsilon$ containing $X$, this sum is calculated solely from $X$. \hfill $\blacksquare$

\bigskip
The derivation just given does not use differential topology concepts such as ``transversality'' or ``regular values in general position''. The mapping degree as defined by an integral over the volume form, goes back to Hadamard and Kronecker. Those intending further to investigate this rich area do well to consult a modern treatment of Brouwer degree such as in \cite{Heinz}, or the monograph of Dinca and Mawhin \cite{D-M}.

 \head 7. Milnor--Hopf Vector Field and Conclusion\endhead

We are not actually using the ``Poincar\'e--Hopf Theorem''. But as in the proof of the renowned P--H\. Theorem, one needs to cobble together a particular vector field that then completes this Complex Axis Theorem.

One may consult the observations in \cite{Spivak} vol. III, p\. 301-302. In fact our special vector field is essentially the same as recommended by Milnor in \cite{TFDV} p\.~40, and expanded upon in the book \cite{Morse Theory} p\. 26.

A consistent way to develop the Milnor--Hopf vector field is by means of the de Medeiros construction. We should begin with an example of an endomorphism of $\CC^{n+1}$ in the form of a non-derogatory, non-defective matrix $L$.  Thus $L$ has {\it linear} elementary divisors (Jordan blocks of size 1), with all eigenvalues $\lambda_j$ {\it distinct}. Consider the diagonal matrix
$$L =
\left[\matrix
0 &&&&&0 \\
\hfil & 1 & \\
\hfil&\hfil & 2 & \hfil \\
\hfil&\hfil&\hfil& \ddots& \hfil \\
\hfil&&&&&n
\endmatrix\right]$$
of size $(n+1) \times (n+1)$. There results the vector field on $\CC^n$, which also has been denoted $\Lambda(L)$,
$$v_L = x_1 \,\frac{\partial}{\partial x_1} + 2x_2 \,\frac{\partial}{\partial x_2}+\cdots+nx_n\,\frac{\partial}{\partial x_n},$$
where the vector field $\hat{v}_L$ is induced on $\CC P^n$ (section of the tangent bundle for complex projective space). The variables $\{x_0,\,\dotsc,\, x_n\}$ are {\it complex}. There is a zero ($\hat{v}_{L}$ becomes radial) at each of the points
$$\align
p_0 &= (1: 0: \cdots : 0) \\
p_1 &= (0:1:\cdots :0)\\
\vdots \\
p_n &= (0: 0: \cdots : 1)
\endalign
  $$
  These points correspond to the totality of the eigen-vectors of $L$. As expected the zeros of $\hat{v}_L$ are simple with index $ = 1$. Indeed, we may work with the typical case of $p_0$ by examining
  $\hat{v}_L$ on the canonical coordinate patch $U_0$ whose ``center'' is $p_0$. One may refer to \cite{Morse Theory} p\. 25-28 or \cite{Hopf} p\. 366-367.

  On this Euclidean patch $U_0 \simeq \CC^n$ we obtain an induced complex vector field
  $$\bold y (z_1,\dotsc, z_n) = (z_1, \,2z_2,\, 3z_3, \dotsc, nz_n)$$
  whose {\it flow} (set of solution curves) emerges as
  $${\bold y}_t(z_1,\dotsc, z_n) = \left( e^t z_1,\, e^{2t} z_2,\, e^{3t} z_3,\dotsc, e^{nt} z_n\right)\; ,  $$
which on a small  poly-disc or real $2n$-ball, gives an orientation-preserving diffeomorphism to its image. Thus the field $y$ has index $ = 1$ at $\vec 0$. The remaining cases treating the other possible ``centers'' $p_1,\,\dotsc,\, p_n$ are handled in the same way with minimal adjustment in notation. Since all these indices are equal and non-zero, Hopf's Lemma shows that the Brouwer degree of $G : \partial W_{\epsilon} \to \SS^{N - 1}$ is also non-zero, so every induced (de Medeiros) field $\hat{v}_{\frak M}$ must have a zero, and every complex square matrix $\frak M$ must have an eigen-vector. This completes a geometric proof of the Complex Axis result, which is also treated in \cite{Derksen} and \cite{de Medeiros}.

  Recapitulating, our proposed square matrix $\frak M$ without eigen-vectors would lead to $\hat{v}_{\frak M}$ with no zeros on $W_{\epsilon}$, which is globally non-degenerate and hence, like $\hat{v}_L$, must have at least one zero after all, contradicting the hypothesis.

\bigskip
\noindent
{\bf Remark}\quad  It would be inconclusive to use certain other matrices besides ``$L$'' as our ``exemplar'', say a matrix that is defective or derogatory. For example, the endomorphism of $\CC^2$ given by
$\left[\matrix
2 & 0 \\
0 & 2
\endmatrix\right]$ leads to a continuum of zeros on $\CC P^1 \simeq \SS^2$, not suitable for counting. The matrix
$\left[\matrix
2 & 1 \\
0 & 2
\endmatrix\right]$ yields an isolated zero, but its index equals 2 (which fact is by now obvious).

\end{document}